\documentclass[10pt,a4paper]{amsart}
\usepackage[latin1]{inputenc}
\usepackage[T1]{fontenc}
\usepackage[english]{babel}

\usepackage{amsmath,amssymb,amsthm}

\theoremstyle{plain}
\newtheorem{theorem}{Theorem}[section]
\newtheorem{corollary}[theorem]{Corollary}
\newtheorem{prop}[theorem]{Proposition}
\newtheorem{lemma}[theorem]{Lemma}
\theoremstyle{definition}
\newtheorem{remark}[theorem]{Remark}
\newtheorem{example}[theorem]{Example}
\newtheorem{definition}[theorem]{Definition}
 \DeclareMathOperator{\re}{Re\,}

\newcommand{\ecc}{\ensuremath{\overline{\mathrm{co}}\,}}
\newcommand{\C}{\mathbb{C}}
\newcommand{\R}{\mathbb{R}}
\newcommand{\K}{\mathbb{K}}
\newcommand{\N}{\mathbb{N}}
\newcommand{\nat}{\mathbb{N}}
\newcommand{\eps}{\varepsilon}
\newcommand{\unifcont}{C_\text{u}}

\renewcommand{\leq}{\leqslant}
\renewcommand{\geq}{\geqslant}

\begin{document}
\title{On the intrinsic and the spatial numerical range}
\thanks{Research partially supported by Spanish MCYT project no.\ BFM2003-01681
and Junta de Andaluc\'{\i}a grant FQM-185}
 \date{January 23th, 2005}

\maketitle

\centerline{\textsc{\large Miguel Mart\'{\i}n,\footnote{Corresponding
author. \emph{E-mail:} \texttt{mmartins@ugr.es}, \emph{Fax:}
\texttt{+34958243272}} } \quad \textsc{\large Javier Mer\'{\i},} \quad
and \quad \textsc{\large Rafael Pay\'{a}} } \vspace{0.3cm}

\begin{center} Departamento de An\'{a}lisis Matem\'{a}tico \\ Facultad de
Ciencias \\ Universidad de Granada \\ 18071 Granada, SPAIN \\
\vspace{0.3cm} \emph{E-mail addresses:} \\
\texttt{mmartins@ugr.es}, \ \texttt{jmeri@ugr.es}, \
\texttt{rpaya@ugr.es}
\end{center}

  \thispagestyle{empty}

\begin{abstract}
For a bounded function $f$ from the unit sphere of a closed
subspace $X$ of a Banach space $Y$, we study when the closed
convex hull of its spatial numerical range $W(f)$ is equal to its
intrinsic numerical range $V(f)$. We show that for every
infinite-dimensional Banach space $X$ there is a superspace $Y$
and a bounded linear operator $T:X\longrightarrow Y$ such that
$\ecc W(T)\neq V(T)$. We also show that, up to renormig, for every
non-reflexive Banach space $Y$, one can find a closed subspace $X$
and a bounded linear operator $T\in L(X,Y)$ such that $\ecc
W(T)\neq V(T)$.

Finally, we introduce a sufficient condition for the closed convex
hull of the spatial numerical range to be equal to the intrinsic
numerical range, which we call the Bishop-Phelps-Bollob\'{a}s
property, and which is weaker than the uniform smoothness and the
finite-dimensionality. We characterize strong subdifferentiability
and uniform smoothness in terms of this property.
\end{abstract}

\section{Introduction}
Given a Banach space $Y$ over $\K$ ($=\R$ or $\C$), we write $B_{Y}$
for the closed unit ball and $S_{Y}$ for the unit sphere of $Y$. The
dual space of $Y$ will be denoted by $Y^{*}$. If $Z$ is another
Banach space, we write $L(Z,Y)$ for the Banach space of all bounded
linear operators from $Z$ into $Y$; if $Z=Y$ we simply write
$L(Y):=L(Y,Y)$ to denote the the Banach algebra of all bounded
linear operators on $Y$. For an element $u\in S_Y$, we write
$$
D(Y,u):=\{y^*\in Y^*\ : \ \|y^*\|=y^*(u)=1\},
$$
the $w^*$-closed and convex set of all \emph{states} of $Y$
relative to $u$. Let us mention two facts, both consequence of the
Hahn-Banach Theorem, which will be relevant to our discussion. On
one hand, we have
\begin{equation}\label{deftau}
\lim_{\alpha \downarrow 0} \frac{\|u+\alpha y\|-1}{\alpha}=
\max\{\re z^*(y) \ : \ z^*\in D(Z,u)\},
\end{equation}
(see \cite[Theorem~V.9.5]{DS} for a proof). On the other hand, if
$X$ is a subspace of $Y$ and $u\in X$, then $D(X,u)$ coincides with
the restriction to $X$ of the elements of $D(Y,u)$.

If $Y$ is a Banach space, by a closed subspace of $Y$ we mean a
Banach space $X$ and an inclusion operator $J:X\longrightarrow Y$
(i.e., $J$ is a linear isometry), and we also say that $Y$ is a
\emph{superspace} of $X$. When no confusion is possible, we omit
$J$, but all the definitions below depend on the way that $X$ is a
subspace of $Y$. Let us fix $X$ and $Y$ as above. We write
$\Pi(X,Y)$ to denote the subset of $S_X\times S_{Y^*}$ given by
$$
\Pi(X,Y):= \left\{(x,y^*)\in S_X\times S_{Y^*} \ : y^*\in
D(Y,Jx)\right\}.
$$
If $X=Y$, we just write $\Pi(Y):=\Pi(Y,Y)$. We denote by
$B(S_X,Y)$ the Banach space of all bounded functions from $S_X$ to
$Y$, endowed with the natural supremum norm, and we write
$\unifcont(S_X,Y)$ for its closed subspace consisting of all
bounded and uniformly continuous functions. For $f\in B(S_X,Y)$ we
can define two different numerical ranges, namely, the
\emph{spatial numerical range} defined as
$$
W(f):=\left\{y^*(f(x))\ : \ (x,y^*)\in \Pi(X,Y)\right\},
$$
and the \emph{intrinsic numerical range} given by
$$
V(f):=\left\{\Phi(f)\ : \ \Phi\in
D\left(B(S_X,Y),J|_{S_X}\right)\right\}.
$$
The name of intrinsic numerical range comes from the fact that if
$f$ belongs to any closed subspace $Z$ of $B(S_X,Y)$, we can
calculate $V(f)$ using only elements in $Z^*$. These numerical
ranges appeared in a paper by L.~Harris \cite{Harris74} for
continuous functions. In the particular case when $X=Y$ and $f$ is
(the restriction to $S_Y$ of) a bounded linear operator, the spatial
numerical range was introduced by F.~Bauer (field of values
subordinate to a norm \cite{Bauer}), extending Toeplitz's numerical
range of matrices \cite{Toe} and, concerning applications, it is
equivalent to Lumer's numerical range \cite{Lumer}. Also in this
case, the intrinsic numerical range appears as the \emph{algebra
numerical range} in the monographs by F.~Bonsall and J.~Duncan
\cite{B-D1,B-D2}; we refer the reader to these books for general
information and background. When $f$ is (the restriction to $S_Y$
of) a uniformly continuous function from $B_Y$ to $Y$ which is
holomorphic on the interior of $B_Y$, both ranges appeared for the
first time in \cite{Harris71}, where some applications are given.

Let us fix a Banach space $Y$ and a closed subspace $X$. For every
$f\in B(S_X,Y)$, $V(f)$ is closed and convex, and we have
\begin{equation}\label{eq-eccWsubsetV}
\ecc W(f)\subseteq V(f),
\end{equation}
where $\ecc$ means closed convex hull. (Indeed, for $x\in S_X$ and
$y^*\in S_{Y^*}$, the mapping $x\otimes y^*$ from $B(S_X,Y)$ to $\K$
defined by $$[x\otimes y^*](g):= y^*(g(x))\qquad \bigl(g\in
B(S_X,Y)\bigr)$$ is an element of $D(B(S_X,Y),J)$.) In the case when
$X=Y$, the inclusion above is known to be an equality whenever $f$
is a uniformly continuous function \cite[Theorem~1]{Harris74} (see
also \cite[\S 9]{B-D1} for bounded linear operators, \cite{Harris71}
for holomorphic functions, and \cite{Rod2004} for a slightly more
general result). On the other hand, the equality $\ecc W(f)=V(f)$
for arbitrary bounded functions cannot be expected in general.
Indeed, this equality holds for every $f\in B(S_Y,Y)$ if and only if
$Y$ is uniformly smooth \cite{Rod2001}. In the general case when $X$
is a proper subspace, two sufficient conditions are given in
\cite[Theorems 2 and 3]{Harris74} for the equality in
Eq.~\eqref{eq-eccWsubsetV}, namely, such a equality holds for all
$f\in \unifcont(S_X,Y)$ if either $X$ is finite-dimensional or $Y$
is uniformly smooth (see definition below). Let us mention that if
$\ecc W(f)= V(f)$ for a bounded function $f\in B(S_X,Y)$, then
$$\max \re V(f)=\sup \re W(f).$$ Therefore, the following formulae,
consequence of Eq.~\eqref{deftau}, will be useful:
\begin{align}
\max \re V(f) &=\lim_{\alpha \downarrow 0} \frac{\|J+\alpha
f\|-1}{\alpha}= \lim_{\alpha \downarrow 0} \ \sup_{x\in
S_X} \frac{\|x + \alpha f(x)\|-1}{\alpha},   \label{igualdad1}\\
\sup \re W(f) &=\sup_{x\in S_X} \lim_{\alpha \downarrow 0} \frac{\|x
+ \alpha f(x)\|-1}{\alpha}. \label{igualdad2}
\end{align}

To state the main results of the paper, let us recall some
definitions and notations.

The norm of a Banach space $Y$ is said to be \emph{smooth} at $u\in
S_Y$ if $D(Y,u)$ reduces to a singleton, and it is said to be
\emph{Fr\'{e}chet-smooth} or \emph{Fr\'{e}chet differentiable} at $u$
whenever there exists
\begin{equation}\label{eq-Frechet}
\lim _{\alpha \rightarrow 0} \frac{\Vert u+\alpha y\Vert -1}{\alpha}
\end{equation}
uniformly for $y\in B_Y$. If this happens for all $u\in S_Y$ we say
that the norm of $Y$ is Fr\'{e}chet differentiable. If, in addition, the
limit in \eqref{eq-Frechet} is also uniform in $u\in S_X$, we say
that the norm of $Y$ is \emph{uniformly Fr\'{e}chet differentiable} at
$S_Y$ or that $Y$ is \emph{uniformly smooth}. A natural succedanea
of Fr\'{e}chet differentiability of the norm when smoothness is not
required is the following notion introduced by D.~Gregory
\cite{Gregory}. The norm of $Y$ is \emph{strongly subdifferentiable}
(\emph{ssd} in short) at $u$ whenever there exists
$$
\lim_{\alpha\downarrow 0} \frac{\|u+\alpha y\|-1}{\alpha}
$$
uniformly for $y\in B_Y$. If this happens for all $u\in S_Y$, we
simply say that the norm of $Y$ is ssd. Thus, the norm of $Y$ is
Fr\'{e}chet differentiable at $u$ if and only if it is strongly
subdifferentiable at $u$, and $Y$ is smooth at $u$. This property
has been fully investigated in \cite{Fran-P}, where we refer the
reader for background. It is shown in \cite[Theorem~1.2]{Fran-P}
that the norm of $Y$ is ssd at $u$ if and only if \emph{$D(Y,u)$ is
strongly exposed by $u$}, i.e., for every $\eps>0$, there exists
$\delta>0$ such that
$$
y^*\in B_{Y^*},\quad\re y^*(u)>1-\delta \qquad \Longrightarrow
\qquad d(y^*,D(Y,u))<\eps.
$$

In this paper we study when the equality in
Eq.~\eqref{eq-eccWsubsetV} holds. The results of the paper can be
divided in two categories.

The first category consists of negative results: we present examples
of pairs of Banach spaces $Y$ and closed subspaces $X$ in which the
equality in Eq.~\eqref{eq-eccWsubsetV} fails, even for elements of
$L(X,Y)$. In section~\ref{sect-finitedimensional} we show that for
every infinite-dimensional Banach space $X$, there is a superspace
$Y$ and an element $T\in L(X,Y)$ such that $\ecc W(T) \neq V(T)$. In
section~\ref{section:FR}, we give concrete examples of Banach spaces
$Y$ for which there is a closed subspace $X$ and an element $T\in
L(X,Y)$ such that $\ecc W(T) \neq V(T)$. Such examples are $c_0$,
$\ell_2\sum_\infty (\ell_2 \oplus_1 \ell_2)$, and, up to renorming,
every non-reflexive Banach space. We will use the following
notation: a Banach space $Y$ is said to have the \emph{FR-property}
if for every closed subspace $X$ and every $T\in L(X,Y)$, the
equality $\ecc W(T)=V(T)$ holds.

The second category is that consisting of positive results. We
introduce in section~\ref{sect-BPB} a sufficient condition for the
FR-property which covers all the previously known examples and may
be interesting by itself. We use the name
``Bishop-Phelps-Bollob\'{a}s property'' for it since it is related
to the quantitative version of the Bishop-Phelps theorem
\cite{BishopPhelps61} given by B.~Bollob\'{a}s \cite{Bollobas}. We
relate this property to the strong subdifferentiability of the norm
and to the uniform smoothness.

\section{When we fix the subspace}\label{sect-finitedimensional}
We recall that, when $X$ is finite-dimensional, for every
superspace $Y$ and every (uniformly) continuous function
$f:S_X\longrightarrow Y$, the equality $\ecc W(f)=V(f)$ holds
\cite[Theorem~2]{Harris74}. The aim of this section is to show
that this fact characterizes the finite-dimensionality, even if we
restrict ourselves to bounded linear operators.

\begin{theorem}\label{theorem-findim}
Let $X$ be an infinite-dimensional Banach space. Then, there are a
superspace $Y$ and an operator $T\in L(X,Y)$ such that $\ecc
W(T)\neq V(T)$
\end{theorem}

We need the following easy lemma.

\begin{lemma}
If $X$ is an infinite-dimensional Banach space, then there exists a
norm-one operator $S \in L(X,c_0)$ which does not attain its norm.
\end{lemma}

\begin{proof}
Since $X$ is infinite dimensional, the Josefson-Nissenzweig theorem
(see \cite[\S XII]{Die}) assures the existence of a sequence
$\{x_n^*\}$ in $S_{X^*}$ $w^*$-converging to $0$. Now, the operator
$S:X \longrightarrow c_0$ defined by
$$[Sx](n)=\dfrac{n}{n+1}x_n^*(x) \qquad \bigl(x\in X,\ n \in \N\bigl),$$
does not attain its norm.
\end{proof}

\begin{proof}[Proof of Theorem~\ref{theorem-findim}]
Let $Y=X \oplus c_0$ endowed with the norm $$\|(x,t)\|=\max
\left\{\|x\|,\ \|Sx\|_\infty+\|t\|_\infty\right\},$$ where $S\in
L(X,c_0)$ is a norm-one operator which does not attain its norm,
and let $J:X\longrightarrow Y$ be the natural inclusion $Jx=(x,0)$
for every $x\in X$. If we define $T\in L(X,Y)$ by $Tx=(0,Sx)$ for
every $x\in X$, it is straightforward to check that
$$\lim_{\alpha \downarrow 0} \ \sup_{x\in
 S_X} \frac{\|x + \alpha T x\|-1}{\alpha} = 1 \qquad \text{and} \qquad \sup_{x\in S_X} \lim_{\alpha
\downarrow 0} \frac{\|x + \alpha T x\|-1}{\alpha}=0.$$ Thus,
Eq.~\eqref{igualdad1} and \eqref{igualdad2} give $V(T) \neq \ecc
W(T)$, as desired.
\end{proof}

\begin{remark}
With a bit more of work, one can show that {\slshape the superspace
$Y$ in the above theorem can be found in such a way that $Y/X$ has
dimension $1$.}\  We divide the proof in two cases, depending on
whether $X$ is reflexive or not.

\noindent \textsc{Case 1:} Suppose $X$ is not reflexive. Then by the
James theorem, there exists $x^* \in S_{X^*}$ which does not attain
its norm. Thus, we can define $Y=X\oplus \K$ endowed with the norm
$$\|(x,t)\|=\max\{\|x\|,\ |x^*(x)|+|t|\} \qquad \bigl(x \in X,\ t\in
\K\bigr),$$ which contains $X$ as the subspace $\{(x,0)\ : \ x\in
X\}$. If we take $T \in L(X,Y)$ defined by $Tx=(0,x^*(x))$ for every
$x\in X$, it is straightforward to show, by using
Eq.~\eqref{igualdad1} and \eqref{igualdad2}, that $\max\re V(T)=1$
and $\sup\re W(T)=0$.

\noindent \textsc{Case 2:} Suppose $X$ is reflexive. By the
Elton-Odell ($1+\varepsilon$)-separation theorem, there are
$\eps_0>0$ and a sequence $\{x_n^*\}_{n\geq 0}$ of elements of
$S_{X^*}$, satisfying $$\|x^*_n - x^*_m\|\geq 1+\eps_0 \qquad
\bigl(n\neq m\bigr)$$ (see \cite[\S XIV]{Die}). Since $X$ is
reflexive, for each $n\in\N$ there exists $x_n\in S_X$ such that
$$|(x_n^*-x_0^*)(x_n)|=\|x_n^*-x_0^*\|\geq 1+\eps_0.$$ Therefore,
\begin{equation}\label{eq-EltonOdell1}
|x^*_0(x_n)|\geq |(x^*_n-x^*_0)(x_n)|-|x_n^*(x_n)|\geq
1+\eps_0-1=\eps_0.
\end{equation}
On the other hand, for each $n\in\N$, we take
$$y_n^*=\dfrac{x_n^*-x^*_0}{\|x_n^*-x^*_0\|}\in S_{X^*}$$ and we observe
that $y_n^*(x_n)=1$ for every $n\in\N$. Since $X\nsupseteq c_0$, it
can be deduced from the proof of the Elton-Odell theorem that
$\{x^*_n\}$ is a basic sequence and so, it converges to zero in the
$weak$ topology by the reflexivity of $X^*$ (see
\cite[Theorem~II.7.2]{Singer}). Using this, and the fact that
$$\|x_n^*-x_0^*\|\geq 1+\eps_0 \qquad \text{and} \qquad
\|x_0^*\|\leq 1,$$ we obtain $$\overline{\lim}\, y_n^*(x)<1 \qquad
\bigl(x\in B_{X}\bigr).$$ This clearly implies that the operator
$S\in L(X,\ell_\infty)$ given by
$$
[Sx](n)=\frac{n}{n+1}y_n^*(x)\qquad \bigl(x \in X,\ n \in
\N\bigr)
$$
does not attain its norm. Now, we take $Y=X\oplus \K$ with the
norm given by
$$
\|(x,t)\|=\max\{\|x\|,\|Sx\|_\infty+|t|\} \qquad \bigl(x\in X,\
 t \in \K\bigr),
$$
we write $J\in L(X,Y)$ for the natural inclusion and, we consider
the operator $T\in L(X,Y)$ defined by $Tx=(0,x_0^*(x))$ for all
$x\in X$. Using Eq.~\eqref{igualdad2} and the fact that $S$ does not
attain its norm, we obtain $\sup \re W(T) = 0$. To compute $\max \re
V(T)$, we observe that
$$
\|J + \alpha T\| \geq \|x_n + \alpha T x_n\|=\|(x_n,\alpha
x_0^*(x_n)\|\geq  \|Sx_n\| + \alpha |x_0^*(x_n)|
$$
so, by using Eq.~\eqref{eq-EltonOdell1} and the fact that
$\|Sx_n\|\longrightarrow 1$, we get
$$
\|J+\alpha T\| \geq 1 + \alpha\eps_0 \qquad \text{for every
$\alpha>0$.}
$$
By just using Eq.~\eqref{igualdad1}, we get $\max \re V(T)\geq
\eps_0$, which finishes the proof. \qed
\end{remark}

\section{When we fix the superspace}\label{section:FR}
As we commented in the introduction, the following result is a
particular case of \cite[Theorems 2 and 3]{Harris74}.

\begin{prop}\label{prop-fd-smooth-FR}
Finite-dimensional spaces and uniformly smooth spaces have the
FR-property.
\end{prop}

In the preceding section we have constructed examples ad hoc of
Banach spaces $Y$ which do not have the FR-property. The aim of this
section is to present some concrete examples of this phenomenon
which will also show that some natural extensions of
Proposition~\ref{prop-fd-smooth-FR} are not possible.

Let us give the first example.

\begin{example}\label{Example-c0noFR}
{\slshape $c_0$ does not have the FR-property.}\ Indeed, let
$Y=c_0\oplus \K^2$ endowed with the norm
$$
\|(x,\lambda,\mu)\|=\max\{\|x\|_\infty, |\lambda| + |\mu|\} \qquad
\bigl(x\in c_0,\ \lambda,\mu\in \K\bigr),
$$
which is isometrically isomorphic to $c_0$. We take a norm-one
functional $x_0^*$ on $c_0$ not attaining its norm, we consider
the closed subspace $X=\{(x,x^*_0(x),0)\ : \ x\in c_0\}$ of $Y$,
and we write $J$ for the natural inclusion of $X$ into $Y$. If we
consider the operator $T:X\longrightarrow Y$ given by
$$
T(x,x_0^*(x),0)= (0,0,x_0^*(x)) \qquad \bigl(x\in c_0\bigr),
$$
by using Eq.~\eqref{igualdad1}, Eq.~\eqref{igualdad2}, and the
fact that $x_0^*$ does not attain its norm, it is easy to verify
that
$$\max\re V(T)=1 \qquad \text{and} \qquad \sup\re W(T)=0,$$
which finish the proof.
\end{example}

Since the norm of $c_0$ is ssd (see \cite[corollary~2.6]{Fran-P},
for instance), the above example shows that
Proposition~\ref{prop-fd-smooth-FR} cannot be extended to the
class of Banach spaces with ssd norm.

On the other hand, using the ideas appearing in the above example,
it is easy to prove the following.

\begin{prop}\label{prop-nonreflexive-nonBPB}
Every non-reflexive Banach space admits an equivalent norm failing
the FR-property.
\end{prop}

\begin{proof}
Let $Z$ be a non-reflexive Banach space. Then, $Z$ is isomorphic to
$Y=V\oplus_\infty(\K\oplus_1\K)$, where $V$ is a $2$-codimensional
closed subspace of $Z$ and, therefore, it is also non-reflexive.
Then, we choose $v_0^*\in S_{V^*}$ which does not attain its norm,
we define the closed subspace
$$X=\left\{(v,v_0^*(v),0) \ : \ v\in V\right\},$$ and we consider $J$ the
natural inclusion of $X$ in $Y$. As in the preceding example, the
operator $T:X\longrightarrow Y$ given by
$$
T(v,v_0^*(v),0)=(0,0,v_0^*(v)) \qquad \bigl(x\in X\bigr)
$$
satisfies
\begin{equation*}
\max\re V(T)=1 \qquad \text{and} \qquad \sup\re W(T)=0.\qedhere
\end{equation*}
\end{proof}

In view of Propositions \ref{prop-fd-smooth-FR} and
\ref{prop-nonreflexive-nonBPB}, one may wonder if reflexivity
implies the FR-property. This is not the case, as the following
example shows.

\begin{example}
{\slshape The superreflexive space $Y=\ell_2\oplus_\infty
\left(\ell_2\oplus_1\ell_2\right)$ does not have the FR
property.}\  \emph{Proof.} First of all, it is straightforward to
show that the norm-one operator $S:\ell_2\longrightarrow\ell_2$
defined by
$$
[Sx](n)=\dfrac{n}{n+1}x(n) \qquad \bigl(x\in \ell_2,\ n\in \N\bigr)
$$
does not attain its norm. Now, we consider the closed subspace
$$
X=\{(x,Sx,0)\ : \ x\in\ell_2\}
$$
with its natural inclusion in $Y$, and we define the operator
$T:X\longrightarrow Y$ by
$$
T(x,Sx,0)=(0, 0, Sx) \qquad \bigl(x\in X\bigr).
$$
The proof will be finished if we show that $\sup\re W(T)=0$ and
$\max\re V(T)\geq1$. For the first equality, given $x\in S_{\ell_2}$
we may find $\alpha_x >0$ such that $(1+\alpha_x)\|Sx\|<1$. Then,
for each $0<\alpha<\alpha_x$ we have
$$
\|(x,Sx,0)+\alpha T(x,Sx,0)\|=\|(x,Sx,\alpha Sx)\|=
\max\{1,(1+\alpha)\|Sx\|\}=1,
$$
and therefore
$$\lim_{\alpha\downarrow 0}\dfrac{\|(x,Sx,0)+\alpha T(x,Sx,0)\|-1}{\alpha}=0.$$
The arbitrariness of $x\in S_{\ell_2}$ gives $\sup \re W(T)=0$. On
the other hand, for each $\alpha>0$, we observe that
$$
\|J+\alpha T\|\geqslant (1+\alpha)\dfrac{n}{n+1} \qquad (n\in \N),
$$
so $\|J+\alpha T\|\geqslant 1+\alpha$, and
$$
\max\re V(T)=\lim_{\alpha\downarrow 0}\dfrac{\|J+\alpha T\|-1}{\alpha} \geqslant 1.
$$
\end{example}

\section[The Bishop-Phelps-Bollob\'{a}s property]{A sufficient condition:
The Bishop-Phelps-Bollob\'{a}s property}\label{sect-BPB} The aim of
this section is to study a sufficient condition for the FR-property
which, actually, covers all the examples given previously. The
motivation for this property is the quantitative version of the
classical Bishop-Phelps' Theorem
\cite{BishopPhelps61,BishopPhelps63} established by B.~Bollob\'{a}s
\cite{Bollobas} (see \cite[\S 16]{B-D2} for the below version).

\begin{theorem}[Bishop-Phelps-Bollob\'{a}s]\label{teorBPB} Let $Y$ be a Banach space and
$\eps>0$. Whenever $y_0\in S_Y$ and $y_0^*\in S_{Y^*}$ satisfy
that $\re y_0^*(y_0)>1-\dfrac{\eps^2}{4}$, there exists
$(y,y^*)\in \Pi(Y)$ such that $$\|y-y_0\|<\eps \qquad \text{and}
\qquad \|y^*-y_0^*\|<\eps.$$
\end{theorem}

This theorem has played an outstanding role in some topics of
geometry of Banach spaces (see \cite{GilesBook,Phelps74,Phelps89},
for instance), specially in the study of ssd norms \cite{Fran-P} or
in the study of spatial numerical range of operators \cite[\S 16 and
\S 17]{B-D2}. Also, the proof of the fact that $\ecc W(f)=V(f)$ for
every $f\in \unifcont(S_Y,Y)$ given in \cite[Theorem~1]{Harris74}
uses the above result. For bounded linear operators, this equality
can be also deduced from \cite[Theorem~8]{Lima93}, a result whose
proof also uses the Bishop-Phelps-Bollob\'{a}s theorem. Motivated by
these facts, we introduce a property which will be sufficient for
the FR-property and it may be of independent interest.

\begin{definition}\label{defBPB}
Let $Y$ be a Banach space and let $X$ be a closed subspace of $Y$.
We say that $(X,Y)$ is a \emph{Bishop-Phelps-Bollob\'{a}s pair}
(\emph{BPB-pair} in short) if for every $\eps>0$ there exists
$\delta >0$ such that whenever $x_0\in S_X$, $y_0^*\in S_{Y^*}$
satisfy $\re y_0^*(x_0)>1-\delta$, there exists $(x,y^*)\in
\Pi(X,Y)$ so that
$$
\|x_0-x\|<\eps \qquad \text{and} \qquad \|y_0^*-y^*\|<\eps.
$$
We say that a Banach space $Y$ has the \emph{BPB property} if for
every closed subspace $X$ of $Y$, $(X,Y)$ is a BPB-pair.
\end{definition}

The next result shows that the BPB property is sufficient for the
FR-property. Actually, it can be proved that the equality in
Eq.~\eqref{eq-eccWsubsetV} holds for uniformly continuous functions.

\begin{theorem}\label{th:BPBsuff}
Let $Y$ be a Banach space and $X$ a closed subspace such that
$(X,Y)$ is a BPB-pair. Then, for every $f\in \unifcont(S_X,Y)$,
the equality $\ecc W(f) = V(f)$ holds.
\end{theorem}

\begin{proof}
Let $J\in L(X,Y)$ be the inclusion map. Let $f\in
\unifcont(S_X,Y)$ and $\Phi \in D(\unifcont(S_X,Y),J)$. By
\cite[Proposition~1]{Harris74}, it suffices to show that
\begin{equation}\label{eq-BPBimpliesFR-1}
\re \Phi(f)\leq\sup\re W(f).
\end{equation}
For each $n\in \nat$,
by using \cite[Lemma~1]{Harris74} we may find $x_n\in S_X$ and
$y_n^*\in S_{Y^*}$ such that
\begin{equation}\label{eq-BPBimpliesFR-2}
\re \Phi(f)\leq\re y_n^*(f(x_n))+1/n
\end{equation}
and $y_n^*(x_n)\longrightarrow 1$. Since $(X,Y)$ is a BPB-pair, it
follows that there exists a sequence
$\bigr\{(\widetilde{x}_n,\widetilde{y}^*_n)\bigl\}_{n\in\nat}
\subseteq \Pi(X,Y)$ such that
$$
\{x_n - \widetilde{x}_n\}_{n\in\nat}\longrightarrow 0 \qquad
\text{and} \qquad \{y_n^* -
\widetilde{y}^*_n\}_{n\in\nat}\longrightarrow 0.
$$
By Eq.~\eqref{eq-BPBimpliesFR-2},
\begin{align*}
\re \Phi(f) &\leq \re \widetilde{y}^*_n(f(\widetilde{x}_n)) + \re
[y_n^*-\widetilde{y}_n^*](f(\widetilde{x}_n)) + \re
y_n^*(f(x_n)-f(\widetilde{x}_n)) + 1/n \\ &\leq \sup \re W(f) +
\|y_n^* - \widetilde{y}_n^*\|\,\|f\|_\infty + \|f(x_n) -
f(\widetilde{x}_n)\| + 1/n
\end{align*}
for all $n\in \nat$. Thus, Eq.~\eqref{eq-BPBimpliesFR-1} follows
from the above and the uniform continuity of $f$.
\end{proof}

It is worth mentioning that the above proof follows the lines of
\cite[Theorem~1]{Harris74}.

\begin{corollary}\label{coro:BPBimpliesFR}
Let $Y$ be a Banach space with the BPB property. Then, $Y$ has the
FR-property.
\end{corollary}

As a consequence of the above corollary and
Theorem~\ref{theorem-findim}, we get the following.

\begin{corollary}
Let $X$ be an infinite-dimensional Banach space. Then, there is a
superspace $Y$ of $X$ such that $(X,Y)$ is not a BPB-pair.
\end{corollary}

The above result implies that not every Banach space $Y$ has the BPB
property. For instance, the examples given in
section~\ref{section:FR} of Banach spaces which do not have the
FR-property also fail the BPB property.

\begin{example}
{\slshape The spaces $c_0$ and $\ell_2 \oplus_\infty (\ell_2
\oplus_1 \ell_2)$ fail the BPB property in their canonical
norms.}\  {\slshape Every non-reflexive Banach space admits an
equivalent norm failing the BPB property.}
\end{example}

On the other hand, if we restrict ourselves to finite-dimensional
subspaces, we get a characterization of the ssd norms.

\begin{prop}\label{nota}
Let $Y$ be a Banach space. Then, the norm of $Y$ is ssd if, and
only if, for every finite-dimensional subspace $X\subseteq Y$, the
pair $(X,Y)$ is BPB.
\end{prop}

\begin{proof}
We suppose first that the norm of $Y$ is ssd. Let $X$ be a
finite-dimensional subspace of $Y$ and let $\eps>0$ be given. Since
the norm of $Y$ is ssd,  \cite[Theorem~1.2]{Fran-P} gives us that
for each $x\in S_X$ there exists $\delta_x>0$ so that
$$
y^*\in S_{Y^*}, \ \ \re y^*(x)>1-\delta_x \qquad \Longrightarrow
\qquad d\left(y^*,D(Y,x)\right)<\eps.
$$
Therefore, if for each $x\in S_X$ we define
$$
A_x=\left\{u\in S_X \ : \ \|u-x\|<\min\left\{\eps,\
\frac{\delta_x}{2}\right\}\right\},
$$
the compactness of $S_X$ assures the existence of $x_1,\dots,x_n \in
S_X$ such that
$$
S_X=\bigcup_{i=1}^{n}A_{x_i}.
$$
Then, $\delta=\min\left\{\frac{\delta_{x_i}}{2} \ : \
i=1,\dots,n\right\}$ satisfies the BPB condition. Indeed, let
$x_0\in S_X$ and $y_0^*\in S_{Y^*}$ be such that
$$
\re y_0^*(x_0)>1-\delta.
$$
Since $x_0\in S_X$, there exists $j\in\{1,\dots,n\}$ so that
$x_0\in A_{x_j}$, that is
$$
\|x_0-x_j\|<\min\left\{\eps,\ \frac{\delta_{x_j}}{2}\right\}.
$$
Therefore, $\re y_0^*(x_j)>1-\delta_{x_j}$ which implies the
existence of $y^*\in D(Y,x_j)$ such that $\|y^*-y_0^*\|<\eps$.

To prove the converse, it is enough to fix $x_0\in S_Y$ and to
show that $x_0$ strongly exposes $D(Y,x_0)$
\cite[Theorem~1.2]{Fran-P}. To do so, let $X$ be the subspace of
$Y$ generated by $x_0$ and, fixed $\eps>0$, let $\delta>0$ be
given by the definition of the BPB for the pair $(X,Y)$ and
$\eps/2$. Suppose now that $y_0^*\in S_{Y^*}$ is such that $\re
y_0^*(x_0)>1-\delta$, then there exists $(x,y^*)\in\Pi(X,Y)$ so
that
$$
\|x-x_0\|<\eps/2 \qquad \text{and} \qquad \|y^*-y_0^*\|<\eps/2.
$$
Since $x \in \text{span}(x_0)$, there exists a modulus-one
$\lambda\in \K$ such that $x=\lambda x_0$. Therefore,
$$
|\lambda-1|=\|\lambda x_0-x_0\|=\|x-x_0\|<\eps/2,
$$
and then,
$$
\lambda y^*\in D(Y,x_0) \qquad \text{and} \qquad \|\lambda y^*
-y_0^*\|\leq \|\lambda y^*-y^*\|+\|y^*-y^*_0\|<\eps/2+\eps/2=\eps,
$$
which finishes the proof.
\end{proof}

Since the norm of any finite-dimensional Banach space is ssd (see
\cite[pp.~48]{Fran-P}), we have the following corollary, which also
implies the first part of Proposition~\ref{prop-fd-smooth-FR}.

\begin{corollary}
Every finite-dimensional Banach space has the BPB property.
\end{corollary}

The other class of spaces with the FR-property given in
Proposition~\ref{prop-fd-smooth-FR} is the one of uniformly smooth
spaces. This result can be also deduced from
Corollary~\ref{coro:BPBimpliesFR}, as the following proposition
shows.

\begin{prop}\label{BPBsuaves}
Every uniformly smooth space has the BPB property.
\end{prop}

\begin{proof}
Let $Y$ be an uniformly smooth space. Then, $Y^*$ is uniformly
convex, so, for every $\eps>0$, we may find $\delta>0$ (the
modulus of convexity of $Y^*$) such that
$$
x^*,\, y^* \in S_{Y^*}, \quad \|x^*+y^*\|>2-\delta \qquad
\Longrightarrow \qquad \|x^*-y^*\|<\eps
$$
(see \cite[Chapter~II]{Beau} for instance). Let $X$ be a subspace of
$Y$, and let $x_0\in S_X$ and $y_0^*\in S_{Y^*}$ be so that $\re
y_0^*(x_0)>1-\delta$. If we consider $y^*\in S_{Y^*}$ such that $\re
y^*(x_0)=1$, we have
$$
\|y^*+y_0^*\|\geq\re (y^*+y_0^*)(x_0)>2-\delta
$$
and, therefore,
$$
\|y^*-y_0^*\|<\eps,
$$
which finishes the proof.
\end{proof}

Observe that, in the above proof, the relation $\eps-\delta$ does
not depend on the subspace. The next result shows that this fact
actually characterizes the uniform smoothness.

\begin{prop}
Let $Y$ be a Banach space with the BPB property in such a way that
the relationship between $\eps$ and $\delta$ in
Definition~\ref{defBPB} does not depend on the subspace $X$. Then,
$Y$ is uniformly smooth.
\end{prop}

\begin{proof}
In view of \cite[Proposition 4.1]{Fran-P}, it is enough to show
that the limit
$$
\lim_{t\downarrow 0}\dfrac{\|u+ty\|-1}{t}=:\tau (u,y)
$$
exists uniformly for $y\in B_Y$ and $u\in S_Y$. Given $\eps>0$, let
$0<\delta<2$ be given by the ``uniform'' BPB property. Now, for
$y\in B_Y$, $u\in S_Y$ and $0<t<\frac{\delta}{2}$, we consider
$$
y_t=\dfrac{u+t y} {\|u+t y\|}\in S_Y \qquad \text{and} \qquad
y_t^*\in D(Y,y_t).
$$
It is immediate to check that $\re y_t^*(u)>1-\delta$ so, if we take
$X=\text{span}(u)$, the BPB property assures the existence of
$(x,z_t^*)\in \Pi(X,Y)$ such that $\|x-u\|<\eps$ and
$\|z_t^*-y_t^*\|<\eps$. Since $x\in\text{span}(u)$, there exists a
modulus-one $\lambda \in\K$ such that $x=\lambda u$. Therefore,
$$
|\lambda-1|=\|\lambda u-u\|=\|x-u\|<\eps,
$$
and then,
$$
\lambda z_t^*\in D(Y,u) \qquad \text{and} \qquad \|\lambda z_t^*
-y_t^*\|\leq \|\lambda
z_t^*-z_t^*\|+\|z_t^*-y_t^*\|<\eps+\eps=2\eps,
$$
Now, the facts
$$
\dfrac{\|u+ty\|-1}{t}=\dfrac{\re y_t^*(u+ty)-1}{t}\leq \re
y_t^*(y)
$$
and $\tau(u,y)\geq \re \lambda z_t^*(y)$ (by Eq.~\eqref{deftau}),
give
$$
0\leq \dfrac{\|u+ty\|-1}{t}-\tau (u,y)\leq
\re y_t^*(y)-\re\lambda z_t^*(y)\leq \|\lambda
z_t^*-y_t^*\|<2\eps,
$$
and the arbitrariness of $\eps>0$ finishes the proof.
\end{proof}

We conclude the paper proving that a pair $(X,Y)$ is a BPB-pair
provided that $X$ is an absolute ideal of $Y$. Let us introduce the
necessary definitions. We refer the reader to \cite[\S~21]{B-D2},
\cite{Paya82}, and references therein for background. A closed
subspace $X$ of a Banach space $Y$ is said to be an \emph{absolute
summand} of $Y$ if there exists another closed subspace $Z$ such
that $Y=X\oplus Z$ and, for every $x\in X$ and $z\in Z$, the norm of
$x+z$ only depends on $\|x\|$ and $\|z\|$. We also say that $Y$ is
an \emph{absolute sum} of $X$ and $Z$. This implies that there
exists an absolute norm on $\R^2$ such that
$$\|x+z\|=|(\|x\|,\|z\|)|_a \qquad \bigl(x\in X,\ z\in Z\bigr).$$
By an \emph{absolute norm} we mean a norm $|\cdot|_a$ on $\R^2$ such
that $|(1,0)|_a=|(0,1)|_a=1$ and $|(a,b)|_a=|(|a|,|b|)|_a$ for every
$a,b\in \R$. Useful results about absolute norms are the following
inequality
$$
\max\{|a|,|b|\}\leq|(a,b)|_a\leq|a|+|b|  \qquad a,b\in \R,
$$
and the fact that absolute norms are nondecreasing and continuous in
each variable. We say that $X$ is an \emph{absolute ideal} of $Y$ if
$X^\perp$ is an absolute summand of $Y^*$, in which case, $Y^*$ can
be identified with $X^*\oplus X^\perp$ with a convenient absolute
sum. It is clear that absolute summands are absolute ideals, but the
converse is not true.

Absolute summands and absolute ideals are generalizations of the
well-known \emph{M-summands}, \emph{L-summands}, \emph{M-ideals},
and the more general class of $L_p$-\emph{summands}
\cite{Behetal,HWW}.

\begin{prop}\label{prop-ideal}
Let $Y$ be a Banach space and let $X$ be a closed subspace. If
$X^*$ is an absolute summand of $Y^*$, then the pair $(X,Y)$ is
BPB. In particular, this occurs when $X$ is an absolute ideal of
$Y$.
\end{prop}

We need the following easy result, which we separate from the
proof of the proposition for the sake of clearness.

\begin{lemma}\label{lemmapair}
Let $E$ be $(\R^2, |\cdot|_a)$ where $|\cdot|_a$ is an absolute
norm. We write $$b_0=\max\{b\geq 0 \ : \ |(1,b)|_a=1\},$$ and we
define
$$
A(\delta)=\{(a,b)\in B_E \ : \ a>1-\delta,\ b\geq b_0\} \qquad
(\delta>0).
$$
Then, for every $\eps>0$ there exists $\delta>0$ such that
$\text{diam}(A(\delta))<\eps.$
\end{lemma}

\begin{proof}
Suppose, for the sake of contradiction, that the result does not
hold. Then, there exists $\eps_0>0$ such that for every $n\in \N$,
$\text{diam}(A(\frac{1}{n}))\geq \eps_0$. So, we may find
$(a_n,b_n)\in A(\frac{1}{n})$ such that $|(a_n,b_n)-(1,b_0)|_a\geq
\frac{\eps_0}{2}$, and thus
\begin{equation}\label{deslema}
\frac{\eps_0}{2}\leq |a_n-1|+|b_n-b_0| \qquad (n\in\N).
\end{equation}
Let $\{(a_{\sigma_n},b_{\sigma_n})\}$  be a convergent subsequence
of $\{(a_n,b_n)\}$, and let $(1,b)\in S_E$ be its limit. By
Eq.~\eqref{deslema} and the fact that
$(a_{\sigma_n},b_{\sigma_n})\in A(\frac{1}{\sigma_n})$, it is
immediate to check that
$$
\frac{\eps_0}{2}\leq |b-b_0| \qquad \text{and} \qquad b\geq b_0.
$$
So, $b$ is strictly bigger than $b_0$, a contradiction.
\end{proof}

\begin{proof}[Proof of Proposition~\ref{prop-ideal}]
There exist a subspace $Z$ of $Y^*$ and an absolute norm
$|\cdot|_a$ on $\R^2$ so that $Y^*=X^*\oplus Z$ and
$$
\|(x^*,z^*)\|=\big|(\|x^*\|,\|z^*\|)\big|_a \qquad \bigl(x^*\in
X^*,\ z^*\in Z\bigr).
$$
For $\eps>0$ fixed, we take $\delta_1>0$ given by the preceding
lemma applied for $\eps/3$, and we define
$$
\delta:=\min\left\{\delta_1, \frac{\eps^2}{36}\right\}.
$$
To finish the proof, for $x_0 \in S_X$ and $y_0^*=(x_0^*,z_0^*)
\in S_{Y^*}$ satisfying
$$
\re y_0^*(x_0)=\re x_0^*(x_0)>1-\delta,
$$
we have to find $(x,y^*)\in \Pi(X,Y)$ so that
$$
\|y^*-y_0^*\|<\eps \qquad \text{and} \qquad \|x-x_0\|<\eps.
$$
To this end, since
$$
\|x_0\|=1=\left\|\frac{x_0^*}{\|x_0^*\|}\right\| \qquad \text{and}
\qquad \re \frac{x_0^*}{\|x_0^*\|}(x_0)\geq \re x_0^*(x_0)
>1-\frac{\eps^2}{36},
$$
we can apply the classical Bishop-Phelps-Bollob\'{a}s Theorem
(\ref{teorBPB}) to $\left(x_0,\frac{x_0^*}{\|x_0^*\|}\right)\in
X\times X^*$ to get $(x,x^*)\in \Pi(X)$ such that
\begin{equation}\label{pairinX}
\quad \left\|x^*-\frac{x_0^*}{\|x_0^*\|}\right\|<\frac{\eps}{3}
\quad \text{and} \quad \|x-x_0\|<\frac{\eps}{3}.
\end{equation}
Now, we distinguish two cases. Suppose first that $\|z_0^*\|\leq
b_0$. Then, we take $y^*:=(x^*,z_0^*)$, which satisfies $\re
y^*(x)=1$ and $\|y^*\|=\big|(1,\|z_0^*\|)\big|_a=1$. Using
Eq.~\eqref{pairinX} and the definition of $\delta$, we get
$$
\|y^*-y_0^*\|=\|x^*-x_0^*\|<\frac{\eps}{3} + \delta < \eps.
$$
So, the pair $(x,y^*)$ satisfies the desired condition.

Suppose otherwise that $\|z_0^*\|>b_0$. In this case, we take
$y^*:=\left(x^*, \frac{b_0}{\|z_0^*\|}z_0^*\right)$, which clearly
satisfies $\re y^*(x)=1=\|y^*\|$. Now, $(1,b_0)$ and $(\|x_0^*\|,
\|z_0^*\|)$ belong to $A(\delta)$ and the diameter of this set is
less than $\eps/3$ by Lemma~\ref{lemmapair}, so we have
$$
\big|\|z_0^*\|-b_0\big|\leq \big|(1,b_0)-(\|x_0^*\|,
\|z_0^*\|)\big|_a< \frac{\eps}{3}
$$
and
\begin{equation*}
\|y^*-y_0^*\|=\big|(\|x^*-x_0^*\|,\|z_0^*\|-b_0)\big|_a\leq
\|x^*-x_0^*\| + \big|\|z_0^*\|-b_0\big|<\eps.\qedhere
\end{equation*}
\end{proof}

By just applying the above proposition and
Theorem~\ref{th:BPBsuff}, we get the following.

\begin{corollary}
Let $Y$ be a Banach space and let $X$ be a closed subspace of $Y$
such that $X^*$ is an absolute summand of $Y^*$ (in particular, if
$X$ is an absolute ideal of $Y$). Then, $$\ecc W(f)=V(f)$$ for
every $f\in \unifcont(S_X,Y)$.
\end{corollary}

An interesting particular case is the case of $M$-embedded and
$L$-embedded spaces. A Banach space $X$ is said to be $M$-embedded
if it is an $M$-ideal of $X^{**}$, and it is $L$-embedded if
$X^{**}=X\oplus_1 Z$ for some closed subspace $Z$ of $X^{**}$.

\begin{corollary}
If $X$ is an $M$-embedded or an $L$-embedded space, then
$(X,X^{**})$ is a BPB-pair.
\end{corollary}

We do not know if the assumption of being $M$-embedded or
$L$-embedded in the above result is superabundant.

\vspace{1cm}

\textbf{Acknowledgment:} The authors would like to express their
gratitude to Pradipta Bandyopadhyay, Catherine Finet, Gilles
Godefroy, and Gin\'{e}s L\'{o}pez for their valuable suggestions and
fruitful conversations about the subject of this paper.

\end{document}